\documentclass[12pt]{amsart}

\usepackage{amssymb,amsmath,amsthm}
\usepackage{amsfonts}
\usepackage{latexsym}
\usepackage{enumerate}
\usepackage{wasysym}
\usepackage{fancyhdr}
\usepackage{epsfig}
\usepackage{breqn}
\usepackage{footnpag}
\usepackage{rotating}
\usepackage{amsfonts}
\usepackage{setspace}
\usepackage{fullpage}
\usepackage{enumitem}
\usepackage{bbold}
\usepackage{comment}
\usepackage{pgf,tikz}
\usepackage{hyperref}
\usepackage{verbatim}
\usepackage{color}
\usepackage{mathrsfs}
\usepackage{soul}
\usetikzlibrary{arrows}
\usepackage{array}

\newtheorem{theorem}{Theorem} 
\newtheorem{lemma}[theorem]{Lemma}
\newtheorem{corollary}[theorem]{Corollary}

\newcommand{\thc}{\operatorname{th}_c}

\title{A sublinear bound on the cop throttling number of a graph}

\author{Anthony Bonato\and Sean English}
\thanks{Ryerson University, Toronto, ON, Canada. \textbf{e-mail:} \{abonato,~sean.english\}$@$ryerson.ca}

\begin{document}
\maketitle

\begin{abstract}
We provide a sublinear bound on the cop throttling number of a connected graph. Related to the graph searching game Cops and Robbers, the cop throttling number, written $\thc(G)$, is given by $\mathrm{th}_c(G)=\min_k\{k+\mathrm{capt}_k(G)\}$, in which $\mathrm{capt}_k(G)$ is the $k$-capture time, or the length of a game of Cops and Robbers with $k$ cops on the graph $G$, assuming both players play optimally.

No general sublinear bound was known on the cop throttling number of a connected graph. Towards a question asked in \cite{CRthrottle}, we prove that $\mathrm{th}_c(G)\leq \frac{(2+o(1))n\sqrt{W(\log(n))}}{\sqrt{\log(n)}},$ where $W=W(x)$ is the Lambert W function.
\end{abstract}

\section{Introduction}

\emph{Graph searching} focuses on the analysis of games and graph processes that model some form of intrusion in a network and efforts to eliminate or contain that intrusion. One of the best known examples of graph searching is the game of {\em Cops and Robbers}, which is a vertex-pursuit game played on graphs that was considered first by Quilliot in \cite{q,q2}, and independently, by Nowakowski and Winkler in \cite{nw}. For a survey of graph searching see \cite{bp,by,fomin}, and see \cite{bonato} for more background on Cops and Robbers.

The game of Cops and Robbers is defined as follows. There are two players, with one player controlling a set of pieces, or \emph{cops}, and the second controlling a single piece, the \emph{robber}. The game is played on a graph, and the pieces are always on vertices of the graph. The game is played over a sequence of discrete time-steps; a \emph{round} of the game is a move by the cops together with a subsequent move by the robber. When it is a player's turn to move, each of their pieces may move to a neighboring vertex or may stay stationary, at the discretion of the player. To begin the game, the cops chose which vertices its pieces initially occupy, then the robber chooses which vertex it will begin at. The cops then move first, followed by the robber; thereafter, the players move on alternate steps. The game is perfect information, with each player seeing all the moves of the other player. The cops win if, after some finite number of rounds, a cop occupies the same vertex as the robber. This is called a \emph{capture}. The robber wins if the robber can evade capture indefinitely. Note that $|V(G)|$-many cops can always win, so the minimum number of cops required to win in a graph $G$ is a well-defined positive integer, named the \emph{cop number} of the graph $G$. We note that the cop number of a graph was first introduced in \cite{af}.

Meyniel's conjecture is a central question concerning the asymptotically largest possible cop number of a connected graph. The conjecture claims that if $G$ is a connected graph of order $n$, then $c(G)=O(\sqrt{n}).$
For more on Meyniel's conjecture, see the survey~\cite{bonato1} and Chapter~3 of \cite{bonato}. The best known asymptotic bound for the cop number of a connected graph $G$ is $c(G)\leq\frac{n}{2^{(1-o(1))\sqrt{\log _2n}}},$ which was proved independently in \cite{fkl,lu,ss}.

For any $k \geq c(G)$ an integer, the \emph{$k$-capture time} of $G$, denoted $\mathrm{capt}_k(G),$ is the minimum number of rounds it takes for $k$ cops to capture the robber on $G$, assuming that all players follow optimal play (that is, the cops attempt to minimize the length of the game, while the robber attempts to maximize it). The $k$-capture time of a graph was first introduced in \cite{boncapt}. In \cite{over}, the authors considered the cases when $k> c(G)$ (that is, there are more cops than the cop number) in the setting of so-called \emph{Overprescribed Cops and Robbers}. A related notion, introduced in \cite{CRthrottle}, is the \emph{cop throttling number} of a graph $G$, denoted by $\thc(G)$ and defined as
\[
\mathrm{th}_c(G) = \min_k\{ k + \mathrm{capt}_k(G)\}.
\]
Here we assume that if $k<c(G)$, then the $k$-capture time is infinite. Hence, the cop throttling number minimizes the number of cops $k$ added to the $k$-capture time, and is connected to positive semidefinite throttling (see \cite{PSDthrottle}). Even in graphs with cop number 1 (such as trees), the cop throttling may be large. In \cite{CRthrottle}, it was shown that $\thc(G) = O(\sqrt{n})$ for several families of graphs $G$ such as trees, unicyclic graphs, and incidence graphs of projective planes.

As $\thc(G) \le \gamma(G) +1,$ where $\gamma(G)$ is the domination number of $G$, it is immediate that $\thc(G) \le n/2 +1,$ where $n$ is the order of $G$. In \cite{CRthrottle}, the authors ask for an asymptotically best possible general upper bound on the cop throttling number. The upper bound in terms of $\gamma(G)$ was the best upper bound known, but $\gamma(G)$ is linear in the number of vertices for many graphs $G$, such as in the cases of $\gamma(P_n)=\gamma(C_n)=\left\lceil\frac{n}3\right\rceil$. Hence, no sublinear bound was known on the cop throttling number of a connected graph. The goal of this note is to provide a partial answer to the question posed in \cite{CRthrottle} by showing that the cop throttling number grows sublinearly with the number of vertices, as presented in Theorem~\ref{main}. Our proof relies on efficiently guarding paths, a recursive decomposition of graphs, and the Moore bound.

All graphs we consider are connected, finite, and undirected. All asymptotics given here will be in terms of $n=|V(G)|$. Given two functions $f=f(n)$ and $g=g(n)$, we will write $f=O(g)$ if there exist some constant $c$ such that $f\leq cg$, we write $f\ll g$ or $f=o(g)$ if $f/g\to 0$ as $n\to\infty$, and we write $g\gg f$ if $f\ll g$. All logarithms will be base $e$. For background on graph theory not discussed here, see \cite{west}.

\section{Main Result}

Given a connected graph $G$, a $u$-$v$-\emph{geodesic} is a shortest path between vertices $u$ and $v$. A geodesic is a path that is a $u$-$v$ geodesic for some choice of $u$ and $v$. Let $W=W(x)$ be the \emph{Lambert W function} or \emph{product-log function}, which is the inverse of $y=xe^x$ ($xe^x$ here will be restricted to the domain $x\geq 0$, on which $xe^x$ is injective, so $W$ is well-defined).

An induced subgraph $H$ of $G$ is \emph{$k$-guardable} if after finitely many moves, $k$ cops can arrange themselves in $H$ so that the robber is immediately captured upon entering $H$. For example, a clique is $1$-guardable. At some round, we will say a $k$-guardable subgraph $H$ is \emph{guarded} if for the rest of the game, some set of cops in $H$ stay in position to immediately capture the robber upon entering $H$.

We need two lemmas on guarding paths.

\begin{lemma}\label{lemma guard path exact}
If $P$ is a geodesic of length $k$, then for any $r\geq 1$, we can place $\left\lceil\frac{k+1}{2r+1}\right\rceil$ cops on $P$ such that $P$ will be guarded in at most $r$ steps. Further, after these $r$ steps, only one cop is necessary to continue guarding $P$.
\end{lemma}

\begin{proof}
Given $P=(v_1,v_2,\dots,v_{k+1})$, for each $0\leq j\leq \left\lceil\frac{k+1}{2r+1}\right\rceil-1$, we place one cop at $v_{r+1+(2r+1)j}$. Note that every vertex on $P$ is within distance $r$ from some cop.
	
Assume that the robber is placed somewhere on the graph. From a classic result of \cite{af}, a single cop can guard $P$. We will imagine for a moment that there is a third player that will control a single cop and take their turn immediately after the robber's turn. This third player will place and move their cop in such a way that $P$ is guarded by this cop for the entirety of the game. Our goal is to move one of our cops onto the third player's cop, in which case our cop can then guard $P$. Note that we will not allow this third player to interfere with the robber, so the presence of the third player does not affect the game between the first two players.

Suppose that on our turn, the third player's cop has guarded $P$ and is currently on $v_{i_1}$. If no cop is at $v_{i_1}$ currently, then in the first move, every cop will move along $P$ towards $v_{i_1}$. After the robber moves, the third player moves their cop to some vertex $v_{i_2}$ on $P$ within distance $1$ of $v_{i_1}$. If no cop is at $v_{i_2}$ currently, then every cop will take a step towards $v_{i_2}$. We continue this process recursively. After $j$ steps, given the current position of the robber, the third player moves their cop to $v_{i_j}$ on the path within distance $1$ of $v_{i_{j-1}}$, guarding $P$. If there is not currently a cop on $v_{i_j}$, then every cop will move along $P$ towards $v_{i_j}$. Since every vertex of $P$ started out at distance at most $r$ from a cop, after at most $r$ steps, some cop must land on the third player's cop, and thus, be in position to guard the path. At this point, only this one cop is necessary to continue guarding $P$.
\end{proof}

It is straightforward to see that Lemma~\ref{lemma guard path exact} is sharp since if only $\left\lceil\frac{k+1}{2r+1}\right\rceil-1$ cops are placed on a path with $k+1$ vertices, there will be a vertex at distance at least $r+1$ from every cop. This level of precision has a negligible effect on the proof of Theorem~\ref{main}, however, so we state an immediate corollary of this result that is weaker but easier to use.

\begin{lemma}\label{lemma guard path}
If $P$ is a geodesic of length $r\ell$ for some integers $r,\ell\geq 1$, then we can place $\ell$ cops on $P$ such that $P$ will be guarded in at most $r$ steps, and after these $r$ steps, only one cop is necessary to continue guarding $P$.
\end{lemma}

\begin{proof}
The proof follows immediately from Lemma~\ref{lemma guard path exact}, and the fact that $\left\lceil\frac{r\ell+1}{2r+1}\right\rceil\leq \ell$ for all $r,\ell\geq 1$.
\end{proof}

We now arrive at our main result, which provides a sublinear bound on the cop throttling number of a graph.

\begin{theorem}\label{main}
If $G$ is a connected graph on $n$ vertices, then
	\[
	\mathrm{th}_c(G)\leq \frac{(2+o(1))n\sqrt{W(\log(n))}}{\sqrt{\log(n)}}.
	\]
\end{theorem}

The following bound on the cop throttling number is slightly worse than the one in Theorem~\ref{main}, but it uses only elementary functions.

\begin{corollary}
		If $G$ is a connected graph on $n$ vertices, then
	\[
	\mathrm{th}_c(G)\leq \frac{(2+o(1))n}{(\log(n))^{1/2-o(1)}}.
	\]
\end{corollary}

\begin{proof}
	We claim that $(\log(x))^{\frac{\log\log\log(x)}{\log\log(x)}}\gg \sqrt{W(\log(x))}$. Since $\log\log\log(x)/\log\log(x)=o(1)$, the result will follow from Theorem~\ref{main}.

If $y=(\log(x))^{\frac{\log\log\log(x)}{\log\log(x)}}$, then $\log(y)=\log\log\log(x)$. Hence, $x=\exp(e^y)$. If $z=\sqrt{W(\log(x))}$, then $z^2=W(\log(x))$, so $z^2e^{z^2}=\log(x)$, and finally $x=\exp (z^2e^{z^2})$. It is evident that $\exp(e^x)\ll \exp(x^2e^{x^2})$, and so we have that $(\log(x))^{\frac{\log\log\log(x)}{\log\log(x)}}\gg \sqrt{W(\log(x))}$.
\end{proof}

\medskip

\begin{proof}[Proof of Theorem~\ref{main}]
Let $\tau=\sqrt{\frac{\log(n)}{W(\log(n))}}$ and $\beta=\tau^{\tau^2}$. Let $G$ be a connected graph on $n$ vertices. First, let us consider the case where $\mathrm{diam}(G)\geq \beta\tau$.
	
We will describe how to place cops on $G$ via a recursive algorithm that decomposes $G$ into paths of length $\beta\tau$, stars, paths of length $\tau^2$, and small connected subgraphs. The paths and stars will be guarded with cops, and then we will show that there are enough free cops close to the small connected subgraphs to quickly catch the robber in any of these connected subgraphs.
	
Let $G_1=G$ and let $P_1$ be a geodesic in $G_1$ of length $\beta\tau$. Place $\beta$ cops along $P_1$ according to Lemma~\ref{lemma guard path} to guarantee $P_1$ can be guarded in $\tau$ steps. Let $G_2$ be the graph induced on $V(G_1)\setminus V(P_1)$. Now, recursively for as long as we can, let $P_i$ be a geodesic in $G_i$ of length $\beta\tau$. Place $\beta$ cops along $P_i$ according to Lemma~\ref{lemma guard path}, and let $G_{i+1}$ be the induced subgraph on $V(G_i)\setminus V(P_i)$. We can continue for, say $\ell_1$ steps, until every component in $G_{\ell_1}$ has diameter less than $\beta\tau$. Note that every vertex in $V(G_{\ell_1})$ is distance at most $\beta\tau$ from some path $P_i$.
	
We describe how to cover any large stars in $G_{\ell_1}$. Recursively, for $i\geq \ell_1$, let $v_i$ be a vertex of degree at least $\tau$ in a component of $G_i$. Place a cop at $v_i$ to guard $N(v_i)$, and let $G_{i+1}$ be the subgraph induced on $V(G_i)\setminus V(N(v_i))$. We can continue this until we reach some $G_{\ell_2}$ with $\Delta(G_{\ell_2})<\tau$.
	
We now will find paths of length $\tau^2$. Recursively for $i\geq \ell_2$, let $P_i$ be a geodesic in $G_i$ of length $\tau^2$. We will place $\tau$ cops on $P_i$ according to Lemma~\ref{lemma guard path} so that $P_i$ can be guarded in at most $\tau$ steps. Let $G_{i+1}$ be the induced subgraph on $V(G_i)\setminus V(P_i)$. We can continue this process until we reach some graph $G_{\ell_3}$, such that every component has diameter less than $\tau^2$. This completes the initial placement of the cops. Note that each cop covered at least $\tau$ vertices on average, so the total number of cops used is at most $n/\tau$. Now we will describe how to move the cops to capture the robber quickly.
	
We will guard each of the paths $P_i$ one at a time in order using the cops that were placed on the paths. It is worth noting that $P_i$ may not be a geodesic in $G$, and so $P_i$ may not initially be guardable by a single cop. Once all the vertices in $V(G)\setminus V(G_i)$ are guarded though, the robber is forced to play on $G_i$, in which $P_i$ is a geodesic, and thus, $1$-guardable.

By Lemma~\ref{lemma guard path}, each path takes at most $\tau$ steps to guard. Since each path is of length at least $\tau^2$, there are at most $n/\tau^2$ paths, so this takes at most $\frac{n}{\tau^2}\cdot \tau=n/\tau$ rounds. Once each path has been guarded, if the robber has not been caught yet, the robber must be in a component of $G_{\ell_3}$.
	
By the Moore bound, since $\Delta(G_{\ell_3})<\tau$ and the diameter of every component of $G_{\ell_3}$ is less than $\tau^2$, each component of $G_{\ell_3}$ has order at most $s$, where
	\begin{align*}
	s&=1+\sum_{i=1}^{\tau^2}\tau(\tau-1)^{i-1}\\
	&=(1+o(1))\tau^{\tau^2}<2\beta-2.
	\end{align*}
Since the domination number of a component with $s$ vertices is at most $s/2$, $\beta-1>s/2$ cops can guard whichever component the robber ends up in. By construction, there is a path of length $\beta\tau$ with $\beta$ cops within distance $\beta\tau$ of every vertex in this component. By Lemma~\ref{lemma guard path}, only one cop need remain on each path to keep them guarded, so the $\beta-1$ other cops on this path can then guard the component containing the robber. This takes at most $2\beta\tau$ more steps and the robber is caught. Note that
	\begin{align*}
	\tau(2\beta\tau)&\ll \tau^{2\tau^2}\\
	&=\left(\sqrt{\frac{\log(n)}{W(\log(n))}}\right)^{\left(\frac{2\log(n)}{W(\log(n))}\right)}\\
	&=\left(\sqrt{\frac{W(\log(n))\exp(W(\log(n)))}{W(\log(n))}}\right)^{\left(\frac{2\log(n)}{W(\log(n))}\right)}\\
	&=\exp\left(\frac{1}2~W(\log(n))\frac{2\log(n)}{W(\log(n))}\right)\\
	&=n,
	\end{align*}
so $2\beta\tau\ll n/\tau$. Hence, the total number of rounds it takes to capture the robber with $n/\tau$ cops is at most $(1+o(1))n/\tau$, completing the proof of this case.
	
If the diameter of $G$ is less than $\beta\tau$, then we proceed identically as in the first case, except we do not need to look for geodesics of length $\beta\tau$, and instead proceed immediately to covering large stars, and then geodesics of length $\tau^2$. We will again arrive at a graph with components that have small maximum degree and small diameter, and therefore, have order at most $s$. We then arbitrarily choose a vertex to place $\beta$ cops, and this uses at most $n/\tau+\beta$ cops in all. We move cops identically to the previous case, guarding all the paths in at most $n/\tau$ rounds, and then the $\beta>s/2$ cops placed arbitrarily can then move to and guard whichever component the robber is in after at most $\beta\tau$ more steps due to the small diameter of the original graph. Since $\beta\tau+\beta\ll \frac{\tau^{2\tau^2}}\tau=n/\tau$, adding these together gives a bound of $(2+o(1))n/\tau$, finishing the proof.
\end{proof}

\end{document}